%% file: main.tex
\documentclass[12pt]{article}
\usepackage{amssymb,amsthm,caption2,vmargin}

\input epsf

\setpapersize{USletter}

\setmarginsrb{1in}{1.5in}{1in}{1in} {0pt}{1mm}{1em}{3em}
\usepackage{amsmath,amsfonts}
\usepackage[mathscr]{eucal}
\usepackage{amssymb,amsmath}
\usepackage{latexsym}
\usepackage{amsfonts}
\usepackage{amscd}
\usepackage{graphicx}
\usepackage{mathptmx}
\newcommand {\del}{\partial}
\newcommand {\dbar}{{\bar{d}}}
\newcommand {\xbar}{{\bar{x}}}

\newcommand {\bt}{\tilde{\beta}}
\newcommand {\dbeta}{\hat{\beta}}
\newcommand {\C}{\mathbb {C}}
\newcommand {\R}{\mathbb {R}}
\newcommand {\D}{\mathbb {D}}
\newcommand {\T}{\mathbb {T}}
\newcommand {\cH}{\mathcal {H}}
\newcommand {\Hilbert}{H}
\newcommand {\M}{\mathcal {M}}

\newcommand {\G}{\mathcal {G}}

\newcommand {\cD}{\mathcal {D}}

\newcommand {\N}{\mathcal {N}}
\newcommand {\Z}{\mathbb {Z}}
\newcommand {\half}{\frac {1}{2}}

\newcommand {\st}{\text{ such that }}

\newcommand {\suml}{\sum\limits}

\newcommand {\supl}{\sup\limits}

\newcommand{\norm}[1]{\left\lVert{#1}\right\rVert}

\newcommand{\diam}{{\rm diam}}

\newcommand{\cen}{{\rm center}}
\newcommand{\dist}{{\rm dist}}
\newcommand{\ball}{{\rm Ball}}
\newcommand{\rad}{{\rm radius}}

\newcommand{\supp}{{\rm supp}}

\newtheorem{theorem}{Theorem}[section]
\newtheorem{lemma}[theorem]{Lemma}
\newtheorem{proposition}[theorem]{Proposition}

\newtheorem{rem}[theorem]{Remark}
\numberwithin{equation}{section}

\begin{document}

\title{Ahlfors-Regular Curves In  Metric Spaces.}
\author{Raanan Schul\\
\tt{schul@math.ucla.edu}}         
\date{}
\maketitle\thispagestyle{empty}
\begin{abstract}
We discuss 1-Ahlfors-regular connected sets in a general metric space and  prove that such  sets are `flat' on most scales and in most locations. Our result is  quantitative, and when  combined with work of I. Hahlomaa,  gives a characterization of 1-Ahlfors regular subsets of 1-Ahlfors-regular curves in  metric spaces.  Our result is  a generalization to the metric space setting of the  Analyst's (Geometric) Traveling Salesman theorems of P. Jones, K. Okikiolu, and G. David and S. Semmes, 
and it can be stated in terms of average Menger curvature.
\end{abstract}

\section{Introduction}
\label{intro}
\input{introduction/introduction.tex}

\section{Acknowledgments}
The author would like to thank Immo Hahlomaa and Pertti Mattila for providing  motivation  to work on this problem. In particular, email correspondence with the former regarding Theorem \ref{17_04_06} before its publication.
The author is grateful to John Garnett for many hours of listening and for his many comments on this essay. 
Finally, the author is thankful for important comments and corrections given by Immo Hahlomaa.
 
\section{Proof of Theorems \ref{new-thm-1} and \ref{new-thm-2}}\label{sect-1}
%
\input{preliminaries/preliminaries.tex}
\input{curvy_arcs/curvy_arcs.tex}
\input{almost_straight_lines/almost_straight_lines.tex}

\section{Modifications for proofs of  Theorems \ref{new-thm-1-var} and \ref{new-thm-2-var}}\label{sect-2}
\input{local_var/local_var.tex}
\section{Appendix}
\input{appendix/parameterization.tex}

\bibliographystyle{alpha}
\bibliography{../../bibliography/bib-file}   

\end{document}

%% file: introduction/introduction.tex
%

We will state our new results in subsection \ref{20_04_06}, but first, we will give some basic definitions and  notation, as well as a description of some known results. 

\subsection{Basic  definitions and notation}
%
%
\subsubsection*{Hausdorff length.}
For a set $K$ we denote by  $\cH^1(K)$ the one dimensional  Hausdorff measure, which we call  \textit{Hausdorff length}.

%

\subsubsection*{$\lesssim$ and $\sim$}
Given two functions $a$ and $b$ into $\R$ we say 
\begin{gather*}
a\lesssim b
\end{gather*} 
with constant $C$, when there exists a constant $C=C_{a,b}$ such that 
\begin{gather*}
a\leq C b.
\end{gather*}
We say that $a\sim b$ if we have $a\lesssim b$ and $b\lesssim a$.
We will allow the constants behind the symbols $\sim and \lesssim$ to depend on the 1-Ahlfors-regularity constant (which will be defined later) and the constant $A$ in the definition of $\G^K$ (see equation \eqref{29_11_05}).
\subsubsection*{
Balls and nets,  multiresolution families}
Let $\M$ be a  metric space  with metric $\dist(\cdot,\cdot)$.    
A ball $B$ is a set 
\begin{gather*}
B=\ball(x,r):=\{y:\dist(x,y)\leq r\}
\end{gather*}
for some $x\in\M$ and some $r>0$.
The set
\begin{gather*}
\{y:\dist(x,y)\leq \lambda r\}
\end{gather*}
will be then be denoted by $\lambda B$.  

We say that $X\subset K$ is an $\epsilon-net$ for $K$ if\\
(i) \quad for all $x_1,x_2\in X$ we have $\dist(x_1,x_2) > \epsilon$\\
(ii) \quad for all $y\in K$ there exists $x\in X$ such that $\dist(x,y)\leq \epsilon$\\
Hence $K\subset \bigcup\limits_{x\in X}\ball(x,\epsilon)$ for an $\epsilon-net$ $X$ for $K$.  

Fix a set $K$. Denote by $X^K_n$ a  sequence of $2^{-n}-nets$ for $K$.
Set
\begin{gather}\label{29_11_05}
\G^K=\{\ball(x,A2^{-n}):x\in X^K_n, n \text{ an integer}\}
\end{gather}
for a constant $A>1$.  
Note that we do not assume in this essay  that $X_n\subset X_{n+1}$. 
We call $\G^K$ a {\it multiresolution family}.
Also note that $\G^K$ depends on $K$.

\begin{rem}\label{nesting}
One of the results we quote (Theorem \ref{thm-1}), for which we use this definition of $\G^K$,  requires the additional property that $X_n\subset X_{n+1}$.  To get this we may construct the sets $X_n$ inductively, however we then require some starting point, which we denote by $n=n_0$. 
For Theorem \ref{thm-1} we also require $n_0$ to be sufficiently negative, namely we need
$2^{-n_0}\geq \diam(K)$.
\end{rem}
\subsubsection*{Lipschitz functions, rectifiable sets, rectifiable curves}
A function $f:\R^k \to \M$  is said to be \textit{C-Lipschitz} if for any $x,y\in \R^k$ such that $x\neq y$, 
\begin{gather*}
{\dist(f(x),f(y)) \over \|x-y\|} \leq C\,.
\end{gather*}
A function $f:\R^k \to \M$  is said to be \textit{Lipschitz} if it is \textit{C-Lipschitz} for some $C>0$.
A set is called \textit{k-rectifiable} if it is contained in a countable union of images 
of  Lipschitz functions $f_j:\R^k \to \M$, except for a set of k-dimensional Hausdorff measure zero.
For more details see \cite{Ma}, where one can also find an excellent discussion of rectifiability in the setting of $\R^d$, part of which carries over to  other metric spaces. 

A set is called a \textit{rectifiable curve} if it is the image of a Lipschitz function defined on $\R$. 
\subsubsection*{Geodesic metric spaces}
A minimizing geodesic is  a map $\tau:I\to M$, where $I\subset \R$ is an interval, and $\tau$ preserves distances.
A metric space  is said to be geodesic if any two points  are the two endpoints of  a minimizing geodesic.
\subsubsection*{Ahlfors-regularity}
Given a set $K\subset\M$  we say that $K$ is {\it k-Ahlfors-regular} if there is a constant $C>0$ so that for all $x\in K$ and $0<r<\diam(K)$ we have
\begin{gather*}
{r^k\over C}\leq \cH^k|_K(\ball(x,r))<Cr^k.
\end{gather*}

We say that a connected set $\Gamma\subset \M$ is a 1-Ahlfors-regular curve with constant $C$ if there is a $C>0$ and a surjective C-Lipschitz function $\gamma:[0,1]\to \Gamma$ such that for any $x\in\Gamma$ and $0<r< diam(\Gamma)$ we have
\begin{gather*}
\cH^1(\gamma^{-1}\ball(x,r))\leq Cr.
\end{gather*}
(In this case we automatically have ${r\over C}\leq \cH^1(\gamma^{-1}\ball(x,r))$.)
A 1-Ahlfors-regular curve is often called an Ahlfors-regular curve.

\subsubsection*{The Jones $\beta$ numbers}
Assume we have a set $K$ lying in $\R^d$.
Consider a ball $B$.
We define the Jones $\beta_\infty$ number as 
\begin{eqnarray*}
\beta_{\infty,K}(B)&=&\frac{1}{\diam(B)} \inf\limits_{L \text{ line}}\ 
	\sup\limits_{x \in K \cap B }\dist(x,L)\\
\vspace{3em}\\
&=&{\text{radius of thinnest cylinder containing }K \cap B \over \diam(B)}.
\end{eqnarray*}
Hence if $\hat{K} \supset K$ then $\beta_{\infty,\hat{K}}(B) \geq \beta_{\infty,K}(B)$.
Note that $\beta_\infty$ is scale independent.
This quantity has $L^p$ variants. Given a locally finite measure $\mu$ and $1\leq p <\infty$, one defines
\begin{gather*}
\beta_{p,\mu}(B) = {1 \over \diam (B)} \inf_{L \text{ line}} 
	\left( \int_{B}  \dist(y,L)^p \, {d\mu(y) \over \mu(B)} \right)^{1/p}.
\end{gather*}
Clearly
\begin{gather}\label{27-12-05}
\beta_{p,\mu}\leq \beta_{\infty,\supp(\mu)}
\end{gather}
when the left hand side is defined.
We define  $\beta_{\infty,\mu}=\beta_{\infty,\supp(\mu)}.$

\subsubsection*{Menger curvature and other useful quantities}
Let $x_1,x_2,x_3\in \M$ be three distinct points.
Take $x'_1,x'_2,x'_3\in \C$ such that $\dist(x_i,x_j)=|x'_i-x'_j|$  for $1\leq i,j\leq 3$.
If $x'_1,x'_2,x'_3$ are collinear then define 
\begin{gather*}
c(x_1,x_2,x_3):=0.
\end{gather*}
Otherwise,
let $R$ be the radius of the circle going through  $x'_1,x'_2,x'_3$.  
In this case define
\begin{gather*}
c(x_1,x_2,x_3):={1\over R}.
\end{gather*}
In any case, $c(\cdot)$ is called the Menger curvature.

For an ordered triple $(x_1,x_2,x_3)\in \M^3$ we define
\begin{gather*}
\del_1(x_1,x_2,x_3):=\dist(x_1,x_2)+\dist(x_2,x_3)-\dist(x_1,x_3).
\end{gather*}
Let  $\{x_1,x_2,x_3\}\subset \M$ be  an unordered triple.  Assume without loss of generality 
$\dist(x_1,x_2)\leq\dist(x_2,x_3)\leq\dist(x_1,x_3)$.  Define 
\begin{gather*}
\del(\{x_1,x_2,x_3\}):=\del_1(x_1,x_2,x_3),
\end{gather*}
or equivalently  
\begin{gather*}
\del(\{x_1,x_2,x_3\})=\min\limits_{\sigma\in S_3}\del_1(x_{\sigma(1)},x_{\sigma(2)},x_{\sigma(3)})\,,
\end{gather*}
where $S_3$ is the permutation group on $\{1,2,3\}$.
Hence we have for all $\{x,y,z\}\subset \M$
\begin{gather*}
\del(\{x,y,z\})\leq \diam\{x,y,z\}
\end{gather*}
as well as
\begin{gather*}
0\leq\del(\{x,y,z\})\leq\del_1(x,y,z)\leq 2\diam\{x,y,z\}
\end{gather*}
where non-negativity follows from the triangle inequality.

\begin{rem}\label{17_11_05}
If 
\begin{gather}\label{14_12_05}
\dist(x,y)\leq\dist(y,z)\leq\dist(x,z)\leq A\cdot\dist(x,y)
\end{gather}
then
\begin{gather*}
c^2(x,y,z)\diam\{x,y,z\}^3\sim\del(\{x,y,z\})
\end{gather*}
with constant depending only on $A$. 
Moreover, in a {\it Euclidean} space,  
\begin{gather}\label{17_11_05-2}
\beta_{\infty,\{x,y,z\}}^2(\ball(x, \diam\{x,y,z\}) \diam\{x,y,z\}\sim \del(\{x,y,z\})
\end{gather}
with constant depending only on $A$.
\end{rem}
See \cite{Ha} for the first part of the above remark.  The second part of the remark follows from  the Pythagorean theorem.

We define $\beta_2(B)$ by 
\begin{gather}
\beta^2_2(B)\rad(B)=\int\int\int_{(B\cap\Gamma)^3}\del(\{x,y,z\})\rad(B)^{-3}d\cH^1(z)d\cH^1(y)d\cH^1(x).
\end{gather} 
Note that $0\leq\beta_2(B)\lesssim 1$ 
(where the constant depends only on the 1-Ahlfors-regularity constant).

\subsection{$\R^d$, Hilbert spaces, metric spaces}
We briefly mention some results.  
For more details and historical background see  \cite{D-Park-City}, \cite{Pa1},  the introduction of \cite{DS}, or the survey \cite{my-TSP-survey}.

\begin{theorem}\cite{J1,Ok,my-thesis-as-paper}\label{thm-2}
Let $\Hilbert$ be $\R^d$ or an infinite dimensional  Hilbert space. For any connected set $\Gamma$ and any $K\subset H$  such that $K\subset\Gamma\subset \Hilbert$  we have
\begin{gather}\label{sum_less_length}
	\sum\limits_{\G^K} \beta_{\infty,\Gamma}^2(B)\diam(B) \lesssim \cH^1(\Gamma).
\end{gather}
\end{theorem}
This was first proven for $\R^d$ with $d=2$ by Jones using complex analysis, and then extended to all $d$ by Okikiolu, using geometric methods.  The constant that comes out of Okikiolu's proof depends exponentially on the dimension $d$,  but in  \cite{my-thesis-as-paper}  it was shown that the constants do not depend on the dimension and moreover, that the theorem holds for  an infinite dimensional Hilbert space.  
The following converse theorem gives a very good reason to care about the left hand side of inequality \eqref{sum_less_length}.
\begin{theorem}\cite{J1,my-thesis-as-paper}\label{thm-1}
Let $\Hilbert$ be $\R^d$ or an infinite dimensional  Hilbert space.
Suppose $A$ in the definition of $\G^K$ is large enough, and assume $\G^K$ satisfies the conditions of Remark \ref{nesting} . 
Given a set $K\subset\Hilbert$, there exists a connected set $\Gamma_0\supset K$ such that the length of $\Gamma_0$ satisfies 
\begin{gather}\label{length_leq_sum}
	\cH^1(\Gamma_0) \lesssim 
	\diam(K) + \sum\limits_{\G^K} \beta_{\infty,K}^2(B)\diam(B).
\end{gather}
\end{theorem}
This theorem was shown by Jones for $\R^d$  (\cite{J1}) and, with some modifications, the proof essentially carries over to the setting of an infinite dimensional Hilbert space (see \cite{my-thesis-as-paper}). 
Theorem \ref{thm-1} also has analogues for  general metric spaces (see \cite{Ha,Ha-2}) and for Heisenberg groups (see \cite{FFP}).  

We especially mention  the following  metric space generalization of Theorem \ref{thm-1} for the category of 1-Ahlfors-regular sets.

\begin{theorem}\cite{Ha-3}\label{17_04_06}
Let $K$ be a 1-Ahlfors-regular set 
in a complete geodesic  metric space $\M$ with metric $\dist(\cdot,\cdot)$.
Assume further that for all $z\in K$ and $R>0$
\begin{gather*}
\int\int\int c^2(x_1,x_2,x_3) d\cH^1|_K(x_3) d\cH^1|_K(x_2) d\cH^1|_K(x_1) \leq C_0 R
\end{gather*}
where the integral on the left hand side is over all triples $x_1,x_2,x_3\in K\cap \ball(z,R)$ such that 
\begin{gather*}
A \cdot\dist(x_i,x_j) \geq \diam\{x_1,x_2,x_3\}.
\end{gather*} 
Then there is a  1-Ahlfors-regular connected set $\Gamma_0\supset K$, whose constant depends only on $C_0$ and  on the 1-Ahlfors-regularity constant of $K$.
\end{theorem}

The proof for this theorem is essentially contained in \cite{Ha-2}. Other results of this type and a relevant counterexample are discussed in the survey \cite{my-TSP-survey}.

Before we go on, let us mention an older  result which is  a special case of a much bigger theorem by  David and Semmes.

\begin{theorem}\cite{DS-sing-int-and-rect-sets}\label{DS-SING}
Let $K\subset\R^d$ be a 1-Ahlfors-regular set and $1\leq q\leq \infty$. 
Then $K$ is contained in a connected 1-Ahlfors-regular set if and only if 
for all $z\in K$ and $0<R<\diam(K)$ 
\begin{gather}\label{22_11_05}
\int_0^R\int_{\ball(z,R)}\beta_{q,\cH^1|_K}(\ball(x,t))^2 d\cH^1|_K(x) {dt\over t} \lesssim R.
\end{gather}
\end{theorem}
\begin{rem}
Note  that  the left hand side of inequality \eqref{22_11_05} can be discretized  as a multiresolution sum  as in the left hand side of inequality \eqref{sum_less_length}.  
\end{rem}
The purpose of  this paper is to prove the converse of Theorem \ref{17_04_06}, and thus to obtain  a metric space analogue of Theorem \ref{DS-SING}.

\subsection{New results}\label{20_04_06}
In Section \ref{sect-1} we show the following.
\begin{theorem}\label{new-thm-1}
Let $\Gamma\subset \M$ be a connected 1-Ahlfors-regular set in a metric space.   Then 
 \begin{gather}
\int\int\int_{\Gamma^3}\del(\{x,y,z\})\diam\{x,y,z\}^{-3}
d\cH^1|_\Gamma(z) d\cH^1|_\Gamma(y) d\cH^1|_\Gamma(x)
\lesssim\cH^1(\Gamma).
\end{gather}
The constant behind the symbol $\lesssim$ depends only on the  1-Ahlfors-regularity constant of $\Gamma$. 
\end{theorem}

It follows from Theorem \ref{new-thm-1} that 
 \begin{gather}
\int\int\int c^2(x,y,z)\lesssim\cH^1(\Gamma)
\end{gather}
where the integral is taken over triples $x,y,z\in \Gamma$ such that 
$\dist(x,y)\leq\dist(y,z)\leq\dist(x,z)\leq A\cdot\dist(x,y)$.
The constant behind the symbol $\lesssim$ depends only on the choice of $A$ (which can be given any value greater then $1$) and the 1-Ahlfors-regularity constant of $\Gamma$. 

On route we show 
\begin{theorem}\label{new-thm-2}
Let $\Gamma\subset \M$ be a connected 1-Ahlfors-regular set in a metric space. Let $K\subset \Gamma$ and let $\hat{\G^K}$ be a multiresolution family as in equation \eqref{29_11_05}.
Then we have   
 \begin{gather}
\sum\limits_{B\in\hat{\G^K}}\int\limits_B\int\limits_B\int\limits_B
	\del(\{x,y,z\})\rad(B)^{-3}
		d\cH^1|_\Gamma(z) d\cH^1|_\Gamma(y) d\cH^1|_\Gamma(x)
			\lesssim\cH^1(\Gamma).
\end{gather}
The constant behind the symbol $\lesssim$ depends only on the  1-Ahlfors-regularity constant of $\Gamma$ and the constant $A$ in the definition of $\hat{\G^K}$. 
\end{theorem}

In Section \ref{sect-2} we use these theorems to prove the following.
\begin{theorem}\label{new-thm-1-var}
Let $\Gamma\subset \M$ be a connected 1-Ahlfors-regular set in a metric space.  
Let $z\in \Gamma$ and $R>0$.
Then 
 \begin{gather}
\int\int\int_{(\Gamma\cap\ball(z,R))^3}\del(\{x,y,z\})\diam\{x,y,z\}^{-3}
	d\cH^1|_\Gamma(z) d\cH^1|_\Gamma(y) d\cH^1|_\Gamma(x)
		\lesssim R.
\end{gather}
The constant behind the symbol $\lesssim$ depends only on the  1-Ahlfors-regularity constant of $\Gamma$. 
\end{theorem}

\begin{theorem}\label{new-thm-2-var}
Let $\Gamma\subset \M$ be a connected 1-Ahlfors-regular set in a metric space. Let $K\subset \Gamma$ and let $\hat{\G^K}$ be a multiresolution family as in equation \eqref{29_11_05}.
Then we have   for every $z\in \Gamma$ and $R>0$
 \begin{gather}
\sum\limits_{B\in\hat{\G^K}\atop B\subset\ball(z,R)}\int\limits_B\int\limits_B\int\limits_B
	\del(\{x,y,z\})\rad(B)^{-3}
		d\cH^1|_\Gamma(z) d\cH^1|_\Gamma(y) d\cH^1|_\Gamma(x)
			\lesssim R.
\end{gather}
The constant behind the symbol $\lesssim$ depends only on the  1-Ahlfors-regularity constant of $\Gamma$  and the constant $A$ in the definition of $\hat{\G^K}$. 
\end{theorem}

%% file: preliminaries/preliminaries.tex
%
%
\subsection{Preliminaries, Notation and Definitions}\label{preliminaries}

Assume $\Gamma\subset \M$ is a connected 1-Ahlfors-regular set.
If $\cH^1(\Gamma)=\infty$ then there is nothing to prove.  Hence we may assume
$\cH^1(\Gamma)<\infty$.  
Since the statements of the theorems are invariant under isometry, 
we may replace $\M$ by   $\ell_\infty(\Gamma)$ without loss of generality. This follows from the Kuratowski embedding (see \cite{Heinonen_embedding_lec}).
Thus we may assume  that $\M$ is complete, and that
$$\diam(\ball(x,r))\sim \rad(\ball(x,r))=r\,.$$


\begin{lemma}\label{closure-length}
Assume $\Gamma\subset \M$ is connected.  Then $\cH^1(\Gamma)=\cH^1(\Gamma^{closure})$.
\end{lemma}

\begin{lemma}\label{finite_length_then_cpt}
Assume $\Gamma\subset \M$ is a closed connected set with $\cH^1(\Gamma)< \infty$.
Then $\Gamma$ is compact.
\end{lemma}
\noindent
Proofs of these lemmas can be found in the appendix of \cite{my-thesis-as-paper} (where they are stated for a Hilbert space, but the proofs work in the category of  a complete Metric space).

\smallskip

We will denote by $\T$ the one dimensional torus $\R / \Z $.

\begin{lemma}\label{parameterization}
Let $\Gamma \subset \M$ be a compact connected  set of finite $\cH^1$ measure.  
Then there is  a Lipschitz function $\gamma:\T \to \M \st Image (\gamma)=\Gamma$ 
and $\norm{\gamma}_{Lip} \leq 32\cH^1 (\Gamma)$.
Further, if $\Gamma$ is 1-Ahlfors-regular,  then 
\begin{gather}\label{8_10}
{R \over C}\leq\cH^1(\gamma^{-1}(\ball(x,R)))\leq CR \quad
\forall x\in \Gamma,\ 0<R\leq\diam(\Gamma).
\end{gather}
i.e. 
$\gamma$ will be witness to the fact that $\Gamma$ is an 1-Ahlfors-regular curve.
Here $C$ is a constant depending only on the 1-Ahlfors-regularity constant of the set $\Gamma$.
\end{lemma}
The proof of this lemma is a modification of  a proof in the appendix of  \cite{my-thesis-as-paper}.
This modification is given in the appendix of this paper.

Fix $\gamma:\T\to \Gamma$ as assured by the above lemma.
We may assume without loss of generality that $\gamma$ is an arc-length parameterization (by re-parameterizing by arc-length and by globally scaling the metric so that the total arc-length is $1$).
This also gives us that $\diam(\Gamma)\leq 1$. We will use this fixed $\gamma$ throughout this essay.

Let $\tau=\gamma|_{[a,b]}$.  We denote by  $\ell(\tau)$ the arc-length of $\tau$. 
We will also use $\ell$ as a measure on $\M$ obtained as the push-forward by $\gamma$  of the Lebesgue measure on $\T$. 
By \eqref{8_10},  for any integrable function $f$, we have that $\int f d\ell \sim\int f d\cH^1|_\Gamma$. 

As $K\subset \Gamma$ in the formulation of theorem \ref{new-thm-2} is fixed, we denote by
$\widehat{\G}=\widehat{\G^K}$.
Clearly 
\begin{gather*}
\int\int\int_{\Gamma^3}\del(\{x_1,x_2,x_3\})\diam\{x_1,x_2,x_3\}^{-3}
								d\cH^1(x_3)d\cH^1(x_2)d\cH^1(x_1)
\lesssim\\	
\indent \sum\limits_{B\in \widehat{\G}}
	\int\int\int_{(B\cap \Gamma)^3}\del(\{x_1,x_2,x_3\})\diam(B)^{-3}
								d\cH^1(x_3)d\cH^1(x_2)d\cH^1(x_1).
\end{gather*}
Hence Theorem \ref{new-thm-2} implies Theorem \ref{new-thm-1}.

To prove Theorem \ref{new-thm-2}  we will show 
\begin{gather}\label{28_09}
\sum\limits_{B\in \widehat{\G}}
	\int\int\int_{(B\cap \gamma(\T))^3}\del(\{x_1,x_2,x_3\})\diam(B)^{-3}
								d\ell(x_3)d\ell(x_2)d\ell(x_1)
\lesssim\ell(\gamma),
\end{gather}
or equivalently,
\begin{gather*}
\sum\limits_{B\in \widehat{\G}}\beta^2_2(B)\diam(B)
\lesssim\cH^1(\Gamma).
\end{gather*}

\begin{rem}\label{gamma-one-to-one}
We may consider the isometric embedding $e$ 
$$\M=\ell_\infty(\Gamma){\buildrel{e}\over\longrightarrow} \ell_\infty(\Gamma)\times\{(0,0)\}\subset \ell_\infty(\Gamma)\times[-1,1]^2$$ and obtain a sequence of maps $\gamma_n:\T\to  \ell_\infty(\Gamma)\times[-1,1]^2$ such that $\gamma_n$ is one-to-one,  $\gamma_n\to\gamma$ uniformly,  
$\norm{\gamma_n}_{Lip} \leq (1+2^{-n})\norm{\gamma}_{Lip} $, and 
$\gamma_n$ gives a 1-Ahlfors-regular curve with constant uniformly comparable to that of $\gamma$.  
This means that we may assume without loss of generality that $\gamma$ in inequality \eqref{28_09}  is one-to-one.
This will be useful for the proof of Lemma \ref{find_large_distanced}. 
\end{rem}

We define
\begin{gather}
\G = \{B \in \widehat{\G} :  \cH^1|_\Gamma(4B)<{1\over 6}\ell(\Gamma)\}.
\end{gather}

We first consider $\widehat{\G} \smallsetminus \G$.
\begin{lemma}\label{15_05_06}
$\sum\limits_{B \in \widehat{\G} \smallsetminus \G}\beta_2^2(B)\diam(B) \lesssim \ell (\gamma)$.
\end{lemma}
\begin{proof}
Set $L=\ell(\gamma)$ 

Consider balls $B\in\widehat{\G}$ with $\cH^1(4B)\geq{L\over 6}$ and $\rad(B)\leq AL$.
There are at most  $C$ such balls at each scale, and at most $C'$ scales. The constants $C,\ C'$ depend only on the Ahlfors regularity constant of $\Gamma$ and the constant $A$. 

Consider now balls $B\in\widehat{\G}$ with $\rad(B)> AL$.  There is at most one ball $B$ of each scale, and
$$\beta_2^2(B)\diam(B) \lesssim \diam(B)^{-3} L^3\diam(B)\sim L \frac{L^2}{\diam(B)^{2}}\,.$$
Summing over all scales we get 
\begin{gather*}
\suml_{B\in \widehat{\G} \smallsetminus \G} \beta_2^2(B)\diam(B) \lesssim L
\end{gather*}
\end{proof}


We need some more notation.  
Let $E\subset\M$ be a closed set such that $\Gamma\cap(\M\smallsetminus E)\neq\emptyset$.
We define
\begin{eqnarray*}
\Lambda(E)&:=&
	\{\tau=\gamma|_{[a,b]}: [a,b]\subset \T; [a,b] 
		\text { a connected component of }\gamma^{-1}(\Gamma \cap E)\}.
\end{eqnarray*}
We will freely use $\tau\in \Lambda(E)$ as both  a parameterization of an arc (given by restriction of $\gamma$), and its image.
In particular, we will denote by $\diam(\tau)$ the diameter of the image of $\tau$.

Let $B\in \G$ be a ball.  For $\tau\in\Lambda(B)$ we denote by $\tau^i$ the extension of $\tau$ to an arc in $\Lambda(2^i B)$.  We set 
\begin{gather}
\Lambda^i(B):=\{\tau^i:\tau\in\Lambda(B)\}.
\end{gather}
We will only use $i\in \{0,1,2\}$.
Let $\tau:[a,b]\to \Gamma$ be a sub-arc of $\gamma$ (and hence an arc-length parameterization).
We define  the quantity $\bt(\tau)$ by
\begin{eqnarray*}
\bt^2(\tau)\diam(\tau):=\ell(\tau)^{-3}
	\int\limits_a^b\int\limits_x^b\int\limits_y^b
		\del_1(\gamma(x),\gamma(y),\gamma(z))dzdydx.
\end{eqnarray*}
(This is how we define the Jones $\beta$ number of an arc).

The constant $\epsilon_2$ below will be set in section \ref{almost-flat-arcs} and will depend on the 1-Ahlfors-regularity constant.
Consider $\tau\in \Lambda^2(B)$.  
We call $\tau$ \textit{almost flat} iff
\begin{eqnarray*}
\bt(\tau)\leq \epsilon_2\beta_2(B).
\end{eqnarray*}

We denote  the collection of \textit{almost flat} arcs  in $\Lambda^2(B)$ by   
\begin{eqnarray*}
S_B:=\{\tau \in \Lambda^2(B): \bt(\tau) \leq \epsilon_2 \beta_2(B)\}.
\end{eqnarray*}

Set:
\begin{eqnarray*}
\G_2&:=&\{B\in \G : 	\Lambda^2(B)\subset S_B 	\}\\
\G_1&:=&\G \smallsetminus\G_2
\end{eqnarray*}

We note  that $B\in \G_1$ implies the existence of an arc $\tau_B\in \Lambda^2(B)$ with 
$\tau_B\notin S_B$.  We will make use of this special (possibly non-unique) arc later on.

We will have Theorem \ref{new-thm-2} if we prove 
\begin{gather}\label{beta_over_A_i}
\sum\limits_{B \in \G_i}\beta^2_2(B)\diam(B) \lesssim \ell (\Gamma) 
\end{gather} 
for $i\in\{1,2\}$. 
We prove inequality \eqref{beta_over_A_i} for $i=1$   in subsection \ref{curvy_arcs}, and for $i=2$  in subsection \ref{almost-flat-arcs}.
%
%
%
%
%
%
%
%
%
%
%
%
%
%
%
%
%
%
%
%
%
%
%
%
%
%
%
%
%
%
%
%

%% file: curvy_arcs/curvy_arcs.tex
\subsection{Non-Flat Arcs}\label{curvy_arcs}
In this subsection we prove inequality \eqref{beta_over_A_i} for $i=1$.

We have $\gamma:\T\to \Gamma$.  
Identify $\T$ with $[0,1]$ for the purpose of defining $\cD^0$ - a dyadic decomposition of $\T$ given by the standard dyadic decomposition of $[0,1]$.
We also define $\cD^1$ - the dyadic decomposition of $\T$ corresponding to the rotation of $\T$ by ${1\over 3}$ of a full rotation, i.e. $x\to (x+ {1\over 3} )\mod 1$, and then using the standard dyadic decomposition of $[0,1]$.  
The reason for these two filtrations is the following remark, which earns this (standard) idea the name
 {\bf `the one third trick'}.
\begin{rem}\label{28_10_05}
Given a (possibly non-dyadic) interval $J\subset\T$ such that $\diam(J)< \frac16$ there exits
an interval $I$ such that  
$I\in\cD^0\cup\cD^1$,
with 
$J\subset I$ and $\diam(I)\leq 6\diam(J)$. 
\end{rem}
We  also define the arcs (mappings) 
$\gamma^0:[0,1]\to\Gamma$ and $\gamma^1:[0,1]\to\Gamma$ 
using the above identifications of $[0,1]$ with $\T$.  They should be thought of as two ways of {\it cutting} $\gamma$ at a point. We define $\gamma^i(x,y,z):=(\gamma^i(x),\gamma^i(y), \gamma^i(z))$.

Let $B\in\G_1$.  Let $\tau=\tau_B\notin S_B$.
Let $I$ be a dyadic interval (assured by remark \ref{28_10_05}) such that  
$\gamma^i(I)\supset \tau$ and $\diam(I)\leq6\ell(\tau)\lesssim\diam(\tau)$, where $i=i(\tau)$ is one of $0$ or $1$. 
Note that the mapping $\tau\to I$ is at most $K_1$-to-$1$ for some constant $K_1$ 
depending only on the 1-Ahlfors-regularity constant of $\Gamma$ and the constant $A$ in equation \eqref{29_11_05}.
Assume that  we have  $i(\tau)=0$. 
For  numbers $r,v\in[0,1]$  we will look at the mapping 
$\psi^{v,r}:[0,1]\to[0,1]$ given by $\psi^{v,r}(t)= v+ r t \mod 1$. 
Note that there are  exactly $2^k$ choices of $\tilde{v}$ and corresponding $\tilde{I}\in\cD^0$ 
(of size $2^{-k}$) with 
$\psi^{\tilde{v},r}(\tilde{I})=[v, v+ 2^{-k} r]$. 

For an interval $I\subset[0,1]$ write $I=[a(I),b(I)]$. 
\begin{rem}
When doing addition  $\mod 1$, we have (by change of variable) for any $I'$ with $\diam(I')=2^{-k}$ 
\begin{eqnarray*}
&&\sum_{I\in\cD^0\atop \diam(I)=2^{-k}}
\diam(I)^{-3}
	\int\limits_{a(I)}^{b(I)}\ \int\limits_x^{b(I)}\ \int\limits_y^{b(I)}
		\del_1\circ\gamma^0(x,y,z)dzdydx \\
&\leq&
\diam(I')^{-3}\int\limits_{v=0}^1\ 
	\int\limits_{r=0}^1\ 
			\int\limits_{y\in v+rI'}
				\del_1\circ\gamma^0(v+ra(I'),y,v+rb(I')) 
					 dy \cdot \diam(I')dr dv
\end{eqnarray*}
giving 
\begin{eqnarray*}
&&\sum_{I\in\cD^0\atop \diam(I)=2^{-k}}
\diam(I)^{-3}
	\int\limits_{a(I)}^{b(I)}\ \int\limits_x^{b(I)}\ \int\limits_y^{b(I)}
		\del_1\circ\gamma^0(x,y,z)dzdydx \\
&\leq&
\sum_{I\in\cD^0\atop \diam(I)=2^{-k}}
\diam(I)^{-3}\int\limits_{v=0}^1\ 
	\int\limits_{r=0}^1\ 
			\int\limits_{y\in v+rI}
				\del_1\circ\gamma^0(v+ra(I),y,v+rb(I)) 
					 dy \cdot \diam(I)dr  \diam(I)dv\,.
\end{eqnarray*}
\end{rem}

Let $I'=[a,b]\in\cD^0$.
Define
\begin{gather}
\del_d(\gamma^0 \psi^{v,r}(I')) := \del_1(		\gamma^0(v+ra),
									\gamma^0(v+r{a+b\over 2}),
									\gamma^0(v+rb)	).
\end{gather} 
\begin{lemma}
Let $I\in\cD^0$.
Let $v,r\in[0,1]$ be chosen such that $\psi^{v,r}(I)=[x,z]\ni y$. Then
\begin{gather*}
\del_1\circ\gamma^0(x,y,z) \leq 
\sum\limits_{I'\in\cD, I'\subset I \atop y\in \psi^{v,r}(I')}
	\del_d(\gamma^0\psi^{v,r}(I')).
\end{gather*}
\end{lemma}
\begin{proof}
This is just the triangle inequality reiterated.
\end{proof}

\begin{lemma}
Let $r,v\in[0,1]$ be fixed.
Then
\begin{gather*}
\sum\limits_{I'\in\cD^0}
	\del_d(\gamma^0\psi^{v,r}(I')) 
\lesssim
\cH^1(\Gamma)
\end{gather*}
\end{lemma}
\begin{proof}
We have that $v+r\{I'\in\cD^0\}$ is a dyadic filtration contained in $\T$. 
The sum in the statement of the lemma is therefore a sum of a telescoping series, whose partial sums are bounded by the arc-length of $\gamma$.
\end{proof}
Now,
\begin{eqnarray*}
&&\sum_{I\in\cD^0}
\diam(I)^{-3}
	\int\limits_{a(I)}^{b(I)}\ \int\limits_x^{b(I)}\ \int\limits_y^{b(I)}
		\del_1\circ\gamma^0(x,y,z)dzdydx \\
&\leq&
\sum_{I\in\cD^0}
\diam(I)^{-3}\int\limits_{v=0}^1\ 
	\int\limits_{r=0}^1\ 
			\int\limits_{y\in v+rI}
				\del_1\circ\gamma^0(v+ra(I),y,v+rb(I)) 
					 dy \cdot \diam(I)dr \cdot \diam(I)dv\\
&\leq&
\sum_{I\in\cD^0}
\diam(I)^{-3}\int\limits_{v=0}^1\ 
	\int\limits_{r=0}^1\ 
		\sum\limits_{I'\in\cD^0\atop I'\subset I}\ 
			\int\limits_{y\in v+rI'}
				\del_d(\gamma^0\psi^{v,r}(I')) 
					\cdot  dy \cdot \diam(I)dr \cdot \diam(I)dv\\
&=&
\sum_{I\in\cD^0}
\diam(I)^{-3}\int\limits_{v=0}^1\ 
	\int\limits_{r=0}^1\ 
		\sum\limits_{I'\in\cD^0\atop I'\subset I}
				\del_d(\gamma^0\psi^{v,r}(I')) 
					\cdot r \cH^1(I') \cdot \diam(I)dr \cdot \diam(I)dv\\
&=&
\int\limits_{v=0}^1\ 
	\int\limits_{r=0}^1\ 
	\sum_{I\in\cD^0}
		{1\over \diam(I)} 
			\sum\limits_{I'\in\cD^0\atop I'\subset I}
				\del_d(\gamma^0\psi^{v,r}(I')) 
					\cdot r\cH^1(I')  dr dv\\
&=&
\int\limits_{v=0}^1\ 
	\int\limits_{r=0}^1\ 
		\sum_{I'\in\cD^0}\ 
				\sum\limits_{I\supset I'}
					{\cH^1(I')\over \diam(I)} 
						\del_d(\gamma^0\psi^{v,r}(I'))
						\cdot r  dr dv\\
&\lesssim&
\int\limits_{v=0}^1\ 
	\int\limits_{r=0}^1\ 
		\sum_{I'\in\cD^0} 
						\del_d(\gamma^0\psi^{v,r}(I'))
						\cdot r  dr dv\\
&\lesssim& \ell(\Gamma).
\end{eqnarray*}

Similarly,
\begin{gather*}
\sum_{I\in\cD^1}
\diam(I)^{-3}
	\int\limits_{a(I)}^{b(I)}\ \int\limits_x^{b(I)}\ \int\limits_y^{b(I)}
		\del_1\circ\gamma^1(x,y,z)dzdydx
\lesssim
\ell(\Gamma) 
\end{gather*}
Hence
\begin{gather*}
\sum\limits_{\tau_B \atop B\in\G_1} \bt^2(\tau_B)\diam(\tau_B)
\lesssim
\ell(\Gamma).
\end{gather*}

\begin{lemma}
We have inequality \eqref{beta_over_A_i} for $i=1$.
\end{lemma}
\begin{proof}
\begin{gather}
\sum\limits_{B \in \G_1}\beta(B)^2\diam(B) \lesssim 
\sum\limits_{B \in \G_1}\bt(\tau_B)^2\diam(\tau_B) \lesssim 
\ell (\Gamma).
\end{gather}
\end{proof}

%% file: almost_straight_lines/almost_straight_lines.tex
%
%
\subsection{Almost Flat Arcs}\label{almost-flat-arcs}
In this subsection we prove inequality \eqref{beta_over_A_i} for $i=2$. 

\smallskip
This subsection will have two parts.  We first show that for every ball $B\in \G_2$ there exist  two special arcs, $\eta_1(B)\in\Lambda^1(B)$ and  $\eta_2(B)\in \Lambda^2(B)$.  These arcs will have properties useful for the second part of  this subsection, where we construct a bounded weight which will in turn give us the desired result. 
\subsubsection*{Part I}
\begin{lemma}\label{find_large_distanced}
Let $B\in\G_2$. 
Let $\xi\in\Lambda^2(B)$.
If for every  arc $\tau_i \in \Lambda^1(B)$ 
we have
\begin{gather}\label{11_10_05-111}
\ell(\tau_i)^{-1}\int_{\tau_i}\dist(\cdot,\xi)\leq\epsilon_4\beta_2^2(B)\diam(B)
\end{gather}
then for every triple of arcs $\tau_i,\tau_j, \tau_k\in \Lambda^1(B)$ we have
\begin{gather*}
\diam(B)^{-3}\int_{\tau_i}\int_{\tau_j}\int_{\tau_k}\del(\{x,y,z\})d\ell(z)d\ell(y)d\ell(x)
\leq 
C_2(\epsilon_2^2+\epsilon_4)\beta_2^2(B)\diam(B) 
\end{gather*}
where $C_2$ is a constant which depends only on the 1-Ahlfors-regularity constant of $\Gamma$
\end{lemma} 
\begin{proof}
Let $(\gamma(x_1),\gamma(x_2),\gamma(x_3))\in\Gamma^3$ be an ordered triple.
Let $S_3$ be the permutation group on $\{1,2,3\}$.
We define for $\sigma\in S_3$
\begin{gather*}
\del_\sigma(\gamma(x_1),\gamma(x_2),\gamma(x_3)):=
\del_1(	\gamma(x_{\sigma(1)}),
		\gamma(x_{\sigma(2)}),
		\gamma(x_{\sigma(3)})	)\,.
\end{gather*}
We will  let $\sigma$ depend on a triple $\xbar=(x_1,x_2,x_3)$ and  we will denote this by $\sigma_\xbar$.

Recall that $\del(\{\cdot\})$ is a continuous function.
We denote by $\D_{\tau,n}$ the collection of $2^n$ points in the domain of $\tau$, evenly spaced according to arc-length.
Let $N_0=N_0(B)$ be chosen large enough so that for all 
$\tau_i,\tau_j,\tau_k\in (\Lambda^1(B)\cup \{\xi\})$ 
(possibly non-different) and $n_1,n_2,n_3\geq N_0$ 
\begin{equation}\label{22_10_05-2}
\begin{aligned}
\diam(B)^{-3}&\int_{\tau_i}\int_{\tau_j}\int_{\tau_k}\del(\{x,y,z\})d\ell(z)d\ell(y)d\ell(x)
\sim\\ 
&2^{-n_1-n_2-n_3}
	\sum\limits_{x\in\D_{\tau_i,n_1}} \sum\limits_{y\in\D_{\tau_j,n_2}}\sum\limits_{z\in\D_{\tau_k,n_3}}
	\del(\{\gamma(x),\gamma(y),\gamma(z)\}), 
\end{aligned}
\end{equation}
and for all $n\geq N_0$
\begin{gather}\label{22_10_05-1}
\ell(\tau_i)^{-1}\int_{\tau_i}\dist(\cdot,\xi)
\sim
2^{-n}\sum\limits_{x\in\D_{\tau_i,n}}\dist(\gamma(x),\xi)
\end{gather}

Let $\tau_1,\tau_2,\tau_3\in \Lambda^1(B)$.
Write
\begin{gather*}
\D_{\tau_1,N_0} = \{O_1,O_2,...\},
\end{gather*}
where 
\begin{gather*}
\dist(\gamma(O_i),\xi) \leq \dist(\gamma(O_{i+1}),\xi).
\end{gather*}
Now let us assume for a moment that $\dist(\tau_1,\xi)>0$.
Let $N_1$ be chosen such that
\begin{gather*}
2^{-N_1} < \dist(\gamma(O_1),\xi).
\end{gather*}
Take $N=\max\{N_1,N_0\}$.
We define a function $f$ with domain $\D_{\tau_1,N_0}$ taking values of probability measures on $\D_{\xi,N}$ as follows.
We go over the $O_i$'s as ordered by $i$.
Let $F_i$ be the set 
\begin{gather*}
F_i=\{x'\in \D_{\xi,N}: \dist(\gamma(x'),\gamma(O_i)) \leq 2\dist(\gamma(O_i),\xi)\},
\end{gather*}
which is non-empty by our choice of $N_1$. 
Define $f(O_1)$ as the uniform  probability measure on $F_1$.
Given $f(O_1),...,f(O_{k-1})$, define $f(O_k)$ as the probability measure  on $F_k$, so that the measure
\begin{gather}\label{27_11_06}
\sum_{i\leq k} f(O_i) |_{F_k}
\end{gather} 
is as close as possible (in $\sup$ norm!) to $2^N k$ times the uniform distribution  on $F_k$
(this is our way of ensuring that  \eqref{27_11_06} is as uniform as possible).
We have for all $x\in \D_{\tau_1,N_0}$ and $x'\in \supp(f(x))$, 
\begin{gather*}
	\dist(\gamma(x),\gamma(x'))\leq  2\dist(\gamma(x),\xi). 
\end{gather*}
We also have for any $x'\in \D_{\xi,N}$
\begin{gather}\label{20_10_05}
2^{-N_0}\sum\limits_{x\in\D_{\tau_1,N_0}} f(x)\{x'\} \leq C 2^{-N}
\end{gather}
where $C$ is a constant which depends only on the 1-Ahlfors-regularity constant of $\Gamma$. 
To see inequality \eqref{20_10_05}, assume the contrary. 
Let $O_k$ be the last element such that  $f(O_k)\{x'\}$ was positive.
Then by construction of $f(O_k)$, we have that for all $x''\in F_k$ 
\begin{gather*}
\sum\limits_{i\leq k} f(O_i)\{x''\} \geq 
\sum\limits_{i\leq k} f(O_i)\{x'\} \geq
C 2^{-N+N_0}.
\end{gather*}
Summing over $F_k$ we get a total mass of 
\begin{gather*}
\sum\limits_{x''\in F_k}\sum\limits_{i\leq k} f(O_i)\{x''\} \geq 
\sharp F_k\cdot C 2^{-N+N_0}\geq
C2^{N_0}\frac{\dist(\gamma(O_k),\xi)}{\ell(\xi)}\,.
\end{gather*}
All this mass, however, came from $O_i$'s such that 
\begin{gather*}
\dist(\gamma(O_i),\gamma(O_k))\leq 
2\dist(\gamma(O_i),\xi) + \diam(F_k) + 2\dist(\gamma(O_k),\xi)\leq 
10\dist(\gamma(O_k),\xi)
\end{gather*}
and so by enlarging $C$ we get a contradiction to 1-Ahlfors-regularity.
This gives inequality \eqref{20_10_05}. 

We similarly define $f$ on $\D_{\tau_2,N_0}$ and $\D_{\tau_3,N_0}$.
Now,\\
\begin{eqnarray*}
&&\diam(B)^{-3}\int_{\tau_1}\int_{\tau_2}\int_{\tau_3}\del(\{x,y,z\})d\ell(z)d\ell(y)d\ell(x)\\
&\sim& 
2^{-N_0-N_0-N_0}
	\sum\limits_{x\in\D_{\tau_1,N_0}}
	 	\sum\limits_{y\in\D_{\tau_2,N_0}}
			\sum\limits_{z\in\D_{\tau_3,N_0}}
	\del(\{\gamma(x),\gamma(y),\gamma(z)\})\\
&\lesssim&
2^{-N_0}(\sum\limits_{x\in\D_{\tau_1,N_0}}\dist(\gamma(x),\xi)  +
		\sum\limits_{y\in\D_{\tau_2,N_0}}\dist(\gamma(y),\xi)  +
		\sum\limits_{z\in\D_{\tau_3,N_0}}\dist(\gamma(z),\xi))\\
&& \hspace{1em}+ 2^{-N_0-N_0-N_0}
	\sum\limits_{x\in\D_{\tau_1,N_0}}
	 	\sum\limits_{y\in\D_{\tau_2,N_0}}
			\sum\limits_{z\in\D_{\tau_3,N_0}}\\
&& \hspace{5em}			
				\sum\limits_{x'\in \supp f(x)} \sum\limits_{y'\in \supp f(y)}  \sum\limits_{z'\in \supp f(z)}
				f(x)\{x'\}\cdot (f(y)\{y'\}\cdot f(z)\{z'\}\cdot
	\del(\{\gamma(x'),\gamma(y'),\gamma(z')\})\\
&\lesssim& 
2^{-N_0}(\sum\limits_{x\in\D_{\tau_1,N_0}}\dist(\gamma(x),\xi)  +
		\sum\limits_{y\in\D_{\tau_2,N_0}}\dist(\gamma(y),\xi)  +
		\sum\limits_{z\in\D_{\tau_3,N_0}}\dist(\gamma(z),\xi))\\
&& \hspace{3em}
+ C^3 2^{-N-N-N}
	\sum\limits_{x'\in\D_{\xi,N}}
	 	\sum\limits_{y'\in\D_{\xi,N}}
			\sum\limits_{z'\in\D_{\xi,N}}
	\del_{\sigma_{(x',y',z')}}(\gamma(x'),\gamma(y'),\gamma(z')).
\end{eqnarray*}

We have yet to specify the function $\sigma$ and have total freedom in choosing 
its values in $S_3$.
Choose $\sigma_{(x',y',z')}$ such that 
$\sigma_{(x',y',z')}(x',y',z')$ has increasing order when ordered by $\xi$.
From  inequalities \eqref{11_10_05-111}, \eqref{22_10_05-2},  and \eqref{22_10_05-1} we now get the lemma.

The case  $\dist(\tau_1,\xi)=0$ can either be assumed not to happen by using remark \ref{gamma-one-to-one} or by computing the above integrals (sums) as limits of the corresponding integrals (sums) in the loops $\gamma_n$ from remark \ref{gamma-one-to-one}.
\end{proof}
Let $\xi_2(B)\in\Lambda^2(B)$ be an arc containing the center of $B$. 
We  upper bound the size of   $\epsilon_2$ and fix $\epsilon_4$ in the proof of the following lemma.
\begin{lemma}\label{11_11_05}
Let $B\in\G_2$. 
We have an arc $\xi_1(B)\in \Lambda^1(B)$  such that 
\begin{gather*}
\dbar(B):=\ell(\xi_1(B))^{-1}\int_{\xi_1(B)}\dist(\cdot,\xi_2(B))\geq\epsilon_4\beta_2^2(B)\diam(B)
\end{gather*}
\end{lemma}
\begin{proof}
If the contrary is true then by reducing $\epsilon_4$ and $\epsilon_2$ we get a contradiction from the previous lemma and Ahlfors-regularity (the latter bounds the number of triples).
\end{proof}
We define $\dbeta(B)$ by  
\begin{gather*}
\dbeta^2(B) \ell(\xi_1(B)):=
\dbar(B)=\ell(\xi_1(B))^{-1}\int_{\xi_1(B)}\dist(\cdot,\xi_2(B))
\end{gather*}
\begin{rem}\label{26_04_06}
\begin{gather*}
1\geq \dbeta(B)\gtrsim \sqrt{\epsilon_4}\beta_2(B),
\end{gather*}
with constant depending only  on the 1-Ahlfors-regularity constant of $\Gamma$.
\end{rem}

\subsubsection*{Part II}
\begin{lemma}
Let $R>0$ be given.
There is a $P_1=P_1(R)$ such that one can write a disjoint union 
\begin{gather*}
\G=\G^1\cup...\cup \G^{P_1}
\end{gather*}
where for each $1\leq p_1\leq P_1$ and $B_1,B_2\in \G^{p_1}$ with $\rad(B_1)=\rad(B_2)$, 
we have 
\begin{gather*}
\dist(B_1,B_2)\geq R\cdot\rad(B_1).
\end{gather*}
\end{lemma}
\begin{proof}
By 1-Ahlfors-regularity we have for each $B_0\in \G$ 
\begin{gather*}
\sharp\{B\in \G: \quad (R+1)\cdot B\cap ( R+1)\cdot B_0\neq \emptyset,\quad \rad(B)=\rad(B_0)\}\leq C_1
\end{gather*}
where $C_1$ is some constant depending only on the 1-Ahlfors-regularity constant and the choice of $A$ and $R$.
We create the desired disjoint union by going over the balls in order. We set $P_1=C_1$.  
By the pigeon-hole principle a ball $B$ can be placed in at least one collection $\G^{p_1}$ such that the result of the lemma will not be contradicted.
\end{proof}

The choice of $R$ will be a consequence of lemma \ref{575757}.
Fix $1\leq p_1\leq P_1(R)$.  
Let $M>0$ be any positive integer.
Consider  $\Delta_M^{p_1}\subset\G_2\cap \G^{p_1}$ defined by
\begin{gather*}
\Delta_M^{p_1}:=\{B\in \G_2\cap \G^{p_1}: 2^{-M} \leq \frac12 \dbeta^2(B)  < 2^{-M+1}\}\,.
\end{gather*} 
Write $\Delta_M^{p_1}=\Delta_M^{p_1,1}\cup....\cup\Delta_M^{p_1,KM}$ where 
\begin{gather*}
\Delta_M^{p_1,p_2}:=\{B\in \Delta_M^{p_1}: \rad(B)=A2^{-nKM+p_2}, n\in \Z\},\quad 1\leq p_2\leq KM. 
\end{gather*}
Fix $M>0$ and  $1\leq p_2\leq KM$  ($K$ will be fixed later).
Fix $\Delta\subset\Delta_M^{p_1,p_2}$ a finite subset.
Take $B\in \Delta$.  

We define $Q(B)\subset (1+4\cdot 2^{-KM})2B$ as follows.
Set
\begin{eqnarray*}
U_{0}&:=& 2B\\
U_{n+1}&:=& 
        U_{n} \cup \bigcup\left\{2B':
                B' \in \Delta,\  
                		2B'\cap U_n \neq \emptyset,\  \rad(B)\geq\rad(B')  \right\} \\
Q(B)&:=&\bigcup\limits_n U_n.
\end{eqnarray*}

\begin{proposition}\label{Q-size}
$Q(B)\subset (1+4\cdot 2^{-KM})2B$
\end{proposition}
We first consider the following lemma.
\begin{lemma}\label{575757}
Consider any metric space.
Assume $R=R>0$ is sufficiently large.  
Let $0<\delta<\frac13$.
Let $\{B_i\}_1^n$ be a sequence of balls in this metric space so that 
$B_k \cap B_{k+1}\neq \emptyset$,
$\rad(B_i)\in\{\delta^k:\  k \textrm{ integer}\}$, 
and 
$\dist(B_{i_1},B_{i_2})\geq R\cdot d$ for balls $B_{i_1},B_{i_2}$ satisfying  $\rad(B_{i_1})=\rad(B_{i_2})=d$. 
Then for any $x\in \cup B_i$,
$$dist(x,\cen(B_{k_0}))\leq  \rad(B_{k_0})(1+ 2\delta +  ... +2^k\delta^{k} +...)\,,$$
where $k_0$ is chosen so that 
$$\rad(B_{k_0})=\max_{i=1}^{n}\rad(B_i)\,.$$  
\end{lemma}
\begin{proof}
This follows by induction on $n$.

\noindent
--For $n=1$ this is clear as $\dist(x,\cen(B_1))\leq \rad(B_1)$.

\noindent
--For $n=N+1$:    
Consider the  sequence $B_1, B_2, ..., B_{k_0-1}$.
Let $k_1$ be so that $\rad(B_{k_1})={\max_{i=1}^{k_0-1}\rad(B_i)}$.
By induction,  for any $y\in \cup_{i=1}^{k_0-1}B_i$
$$\dist(y,\cen(B_{k_1}))\leq  \rad(B_{k_1})(1+ 2\delta +  ... +2^k\delta^{k} +...)\,.$$ 
Hence, if $R$ is large enough, $\rad(B_{k_1})\neq \rad(B_{k_0})$, and hence
$\rad(B_{k_1})< \rad(B_{k_0})$ which gives
$$\dist(y,\cen(B_{k_1}))\leq  \delta\rad(B_{k_0})(1+ 2\delta +  ... +2^k\delta^{k} +...)\,.$$ 
Similarly for the sequence  $B_{k_0+1}, B_{k_0+2}, ..., B_n$.
We conclude 
\begin{eqnarray*}
\dist(x,\cen(B_{k_0}))
&\leq&  
\rad(B_{k_0})+ 2\delta\rad(B_{k_0})(1+ 2\delta +  ... +2^k\delta^{k} +...)\\
&=&  
\rad(B_{k_0})(1+ 2\delta(1+ 2\delta +  ... +2^k\delta^{k} +...))\\
&=&  
\rad(B_{k_0})(1+ 2\delta +  ... +2^k\delta^{k} +...)\,.
\end{eqnarray*}
This concludes the induction.
\end{proof}

We now prove Proposition \ref{Q-size}. 
\begin{proof}
Recall that the number of balls in $\Delta$ is finite.
We denote by $\delta=2^{-MK}$.
If $K>2$ we have $\delta<\frac13$.
Let $x\in Q$.  Then there is a sequnce of balls $\{B_i\}_1^n$ such that 
$x\in B_n$,
$\frac12B_i\in\Delta$, 
$B_k \cap B_{k+1}\neq \emptyset$ with
$\rad(B_i)\leq \rad(2B)$ and $B_1=2B$.
Using lemma \ref{575757} we get the desired estimate as
$$(1+ 2\delta +  ... +2^k\delta^{k} +...)\leq 1 + 4\cdot2^{-KM}\,.$$
\end{proof}

The family $\{Q(B):B\in \Delta\}$ has the property that if $Q_1$ and $Q_2$ are in it, then if $Q_1\cap Q_2 \neq \emptyset$ we have $Q_1\subset Q_2$ or $Q_2\subset Q_1$.

We write  
\begin{gather}\label{05_11_05-7} 
Q=(\bigcup\limits_i Q^i) \cup R_Q
\end{gather}
where $Q^i$ is maximal such that $Q^i=Q(B^i),\ B^i\in \Delta$ and $Q^i\subsetneq Q$.
We choose $R_Q$ so that all the  unions in equation  \eqref{05_11_05-7} are disjoint.

Let $B\in \G$ be a ball.  for $\tau\in\Lambda(B)$ we denote by $\tau^Q$ the extension of $\tau$ to an arc in $\Lambda({Q(B)})$.  We set 
\begin{gather}
\Lambda^Q(B):=\{\tau^Q:\tau\in\Lambda(B)\}.
\end{gather}
\begin{rem}\label{28_10_05-990}
We have (using regularity) that if $B\in \G_2$ then for all $\tau\in\Lambda^Q(B)$
\begin{gather*}
\bt(\tau)\lesssim\epsilon_2\beta_2(B)\lesssim 
\epsilon_2 \sqrt{\epsilon_4^{-1}}\dbeta(B).
\end{gather*}
\end{rem}
We also denote by $\xi_2(Q)$ a connected component of $\xi_2(B)\cap Q$ which contains the center of $B$.
We will denote by $J_1(Q)$ and $J_2(Q)$ the index sets 
\begin{gather*}
J_1(Q)=\{i:\quad Q^i \cap \xi_1(B)\neq\emptyset\}\\
J_2(Q)=\{i:\quad Q^i \cap \xi_2(Q)\neq\emptyset\}.
\end{gather*}
\begin{rem}\label{11_11_05-3}
By enlarging $K$ if necessary, if $x\in \xi_1(B)$ such that  
$\dist(x,\xi_2(Q))\geq {1\over 4} 2^{-M}\diam(2B)$ and
$x\in Q^j$, then   
$j\in J_1\smallsetminus J_2$.
\end{rem}
\begin{proposition}\label{lots_of_length}
Let $B\in\Delta$ and $Q=Q(B)$.  
Then 
\begin{gather*}
\ell(R_Q)+
\sum\limits_{j} \diam(Q^j)\geq
\ell(R_Q)+
\sum\limits_{j\in J_1 \cup J_2} \diam(Q^j)
	\geq (1+ c''\dbeta(B))\diam(Q)
\end{gather*}
for some constant $c''>0$ depending only on the 1-Ahlfors-regularity of $\Gamma$.
\end{proposition}

Before we can prove this proposition we need two lemmas.

\begin{lemma}
There is a constant $c>0$,
independent of $\epsilon_2$, $K$,  and $M$, so that
for any $Q=Q(B)$ and $\xi_1=\xi_1(B)$, $\xi_2=\xi_2(Q)$  we have
\begin{gather*}
\ell(R_Q\cap \xi_1)+
\sum\limits_{j\in J_1\smallsetminus J_2} \diam(Q^j)
\geq c2^{-M\over 2}\diam(Q)\,.
\end{gather*}
\end{lemma}
\begin{proof}
Let $\dbar=\dbar(B)$.
Assume for a moment
\begin{gather}\label{5-12-06}
\ell(\{x\in \xi_1:\dist(x,\xi_2(B))\geq {\dbar\over 2}\}) \leq \dbeta(B) \ell(\xi_1).
\end{gather}
Then 
\begin{gather*}
(1-\dbeta(B))\ell(\xi_1) {\dbar\over 2} +  \dbeta(B)\ell(\xi_1)d_\infty  
	\geq \dbar\ell(\xi_1),
\end{gather*}
where 
\begin{gather*}
d_\infty=\sup\limits_{x\in\xi_1}\dist(x,\xi_2(B)).
\end{gather*}
Hence
\begin{gather*}
d_\infty \dbeta(B) \geq 
	\dbar - (1-\dbeta(B)){\dbar\over 2}= {\dbar\over 2} + \dbeta(B){\dbar\over 2}
\end{gather*}
or
\begin{gather*}
d_\infty  \geq 
	 \dbeta(B)^{-1} {\dbar\over 2} + {\dbar\over 2}
\end{gather*}
and hence (since $\xi_1$ is connected and we are assuming \eqref{5-12-06}), the diameter of the largest-diameter connected component of 
\begin{gather*}
\{x\in \xi_1:\dist(x,\xi_2(B))\geq {\dbar\over 2}\}
\end{gather*}
is at least  
\begin{gather*}
\dbeta(B)^{-1} {\dbar\over 2} 
		= \half\dbeta(B) \ell(\xi_1).
\end{gather*}
Either way (with or without assumption \eqref{5-12-06}) we have
\begin{gather*}
\ell(\{x\in \xi_1:\dist(x,\xi_2(Q))\geq {\dbar\over 2}\}) \geq 
\ell(\{x\in \xi_1:\dist(x,\xi_2(B))\geq {\dbar\over 2}\}) \geq \half \dbeta(B) \ell(\xi_1)
\end{gather*}
where the first inequality follows from $\xi_2(Q)\subset \xi_2(B)$.
By remark \ref{11_11_05-3} and the definitions of $\dbeta$ and $\xi_1$,
\begin{eqnarray*}
\ell(R_Q\cap \xi_1)+
\sum\limits_{j\in J_1\smallsetminus J_2} \diam(Q^j)
&\gtrsim&
\ell(\{x\in \xi_1:\dist(x,\xi_2(Q))\geq {\dbar\over 2}\})	\\
&\geq&  
\half \dbeta(B) \ell(\xi_1)
\gtrsim \dbeta(B) \diam(Q)
\gtrsim 2^{-M\over 2} \diam(Q).
\end{eqnarray*}
An important thing to note is that all the similarity constants are independent of $\epsilon_2$, $K$,  and $M$ since these are rough lower bounds.
This gives
\begin{gather*}
\ell(R_Q\cap \xi_1)+
\sum\limits_{j\in J_1\smallsetminus J_2} \diam(Q^j)
\geq c2^{-M\over 2}\diam(Q)
\end{gather*}
with $c$ independent of $\epsilon_2$, $K$,  and $M$.
\end{proof}

\begin{lemma}
There is a constant $\epsilon_3>0$ (independent of $M$), which we can make arbitrarily small
by reducing $\epsilon_2$ and increasing $K$,   such that 
for any $Q=Q(B)$ and  $\xi_2=\xi_2(Q)$  we have
\begin{gather*}
\ell(R_Q\cap \xi_2)+
\sum\limits_{j\in J_2} \diam(Q^j)
\geq (1-\epsilon_3 2^{-M\over 2})\diam(Q)\,.
\end{gather*}
\end{lemma}
\begin{proof}
Throughout the proof we assume $\epsilon_2$ is sufficiently small.
Recall that by construction we have 
$$\diam(Q^j)\leq  (1+4\cdot 2^{-KM})2^{-KM}\diam(Q)$$ and
$$\diam(Q)\leq(1+4\cdot2^{-MK})\diam(2B)\,.$$

Let $\xi_{2.0}\in\Lambda(2B)$ be a subarc of $\xi_2(Q)$ containing the center of $B$.
Let $O$ be the center of $B$, and $O_1,O_2$ the entry and exit points of $\xi_{2.0}$ from $2B$. 
Assume without loss of generality that $O_1 < O < O_2$ as ordered by $\xi_2$. 
Consider a ball $\ball(O_1,r)$, with $r\leq \rad(B)$.
Let $O_1^r\in\xi_2$ be the (unique) point s.t. $\dist(O_1^r,O_1)=r$, $O_1<O_1^r<O$,  and any other such point $X$ satisfies $X<O_1^r$.
Symmetrically,   
let $O_2^r\in\xi_2$ be the (unique) point s.t. $\dist(O_2^r,O_2)=r$, $O<O_2^r<O_2$,  and any other such point $X$ satisfies $X>O_2^r$.  

The constants $r_0$ and $C_{r_0}$ will be fixed below, independently of $\epsilon_2$ and $M$.
Suppose for a moment that there is no pair $r_1,r_2\in[0,r_0\diam(Q)]$ such that 
\begin{gather}\label{08_10-11111}
\partial_1(O_1^{r_1},O,O_2^{r_2})<C_{r_0} \epsilon_2\sqrt{\epsilon_4^{-1}} 2^{-M\over 2}\diam(Q)\,.
\end{gather}
Then
\begin{eqnarray*}
\epsilon_2^2 \beta^2(B)\diam(B) &\gtrsim&
      \bt^2(\xi_2)\diam(\xi_2)\\
&\gtrsim&  
\diam(B)^{-3}(C_{r_0}\epsilon_2\sqrt{\epsilon_4^{-1}} 2^{-M\over 2}\diam(Q))^2 
	\cdot r_0\diam(Q)\cdot r_0\diam(Q)\\
&\sim&
C^2_{r_0}r_0^2\epsilon_2^2 \epsilon_4^{-1} 2^{-M}\diam(Q)\\
&\gtrsim& 
C^2_{r_0}r_0^2\epsilon_2^2\beta^2(B)\diam(B)\,.
\end{eqnarray*}
Thus by setting $C_{r_0}$ large enough with respect to $r_0$ we get a contradiction.  
So we let $r_1,r_2\in[0,r_0\diam(Q)]$ be a pair such that \eqref{08_10-11111} holds.
This implies  	
\begin{eqnarray*}
\dist(O_1^{r_1},O_2^{r_2})
&\geq& 
\dist(O_1^{r_1},O) + \dist(O_2^{r_2},O) - C\epsilon_2 \sqrt{\epsilon_4^{-1}} 2^{-M\over 2}\diam(Q)\\
&\geq& 
\rad(2B)-r_1 + \rad(2B)-r_2 - C\epsilon_2 \sqrt{\epsilon_4^{-1}} 2^{-M\over 2}\diam(Q)\\
&\geq& 
\diam(Q)- r_1 -r_2 - 4\cdot2^{-MK}\diam(Q) -C\epsilon_2 \sqrt{\epsilon_4^{-1}} 2^{-M\over 2}\diam(Q)\,.
\end{eqnarray*}

If  $r_1\leq 2\cdot 2^{-MK}\diam(Q)$ define $\xi_{2.1}=\emptyset$.
If  $r_2\leq 2\cdot 2^{-MK}\diam(Q)$ define $\xi_{2.2}=\emptyset$.
Otherwise, we define $\xi_{2.1}$ or $\xi_{2.2}$ as follows.

For points $X,Y\in \xi_2$, we will denote by $X\leadsto Y$ the subarc of $\xi_2$ connecting $X$ and $Y$.  Assume $r_i> 2\cdot 2^{-MK}\diam(Q)$. Let 
$$B_i=\left(1-2\cdot2^{-MK} \frac {\diam(Q)}{r_i}\right)\ball(O_i,r_i)\,.$$
By the definition of $O_i^{r_i}$  we have that 
$$\dist\left( O_1^{r_1}\leadsto O_2^{r_2},\  B_i  \right)\geq
	2\cdot 2^{-MK}\diam(Q)\,.$$
By reducing $r_0$ the balls $\ball(O_1,r_1)$ and $\ball(O_2,r_2)$ have distance at least $2\cdot2^{-MK}\diam(Q)$ from each other.
Define $\xi_{2.i}$ to be the largest-diameter (connected) subarc of $(O_i\leadsto O_i^{r_i})\cap B_i$.  

In either case, $\diam(\xi_{2.i})\geq r_i-2\cdot 2^{-MK}\diam(Q)$.

Denote by $\xi_{2.3}$ the subarc $O_1^{r_1}\leadsto O_2^{r_2}$.  	By the above we have that no $Q^j$ intersects 2 of these subarcs, and that  (by increasing $K$ for the last inequality)
\begin{eqnarray*}
&&\diam(\xi_{2.1})+\diam(\xi_{2.2})+\diam(\xi_{2.3})\\
&\geq& 
r_1 -2\cdot 2^{-MK}\diam(Q) +  r_2 -2\cdot 2^{-MK}\diam(Q)  \\
&& +\ \diam(Q)- r_1 -r_2 - 4\cdot2^{-MK}\diam(Q)- C\epsilon_2 \sqrt{\epsilon_4^{-1}} 2^{-M\over 2}\diam(Q)\\
&\geq&
\diam(Q)- 8\cdot 2^{-MK}\diam(Q)  - C\epsilon_2 \sqrt{\epsilon_4^{-1}} 2^{-M\over 2}\diam(Q)\\
&\geq&
(1-\epsilon_3 2^{-M\over 2})\diam(Q)\,. 
\end{eqnarray*}
Furthermore, since $\diam(\xi_{2.i}\cap Q^j)\leq \diam(Q^j)$, we have
\begin{gather*}
\ell(R_Q\cap \xi_2)+
\sum\limits_{j\in J_2} \diam(Q^j)
\geq (1-\epsilon_3 2^{-M\over 2})\diam(Q)\,.
\end{gather*}
\end{proof}

We now get Proposition \ref{lots_of_length}:
\begin{proof}
\begin{eqnarray*}
\ell(R_Q)+
\sum\limits_{j} \diam(Q^j) 
&\geq&
\ell(R_Q\cap \xi_2)+
\sum\limits_{j\in J_2} \diam(Q^j) +
\ell(R_Q\cap \xi_1)+
\sum\limits_{j\in J_1\smallsetminus J_2} \diam(Q^j)\\
&\geq& 
(1-\epsilon_3 2^{-M\over 2})\diam(Q) + c2^{-M\over 2}\diam(Q)\,.
\end{eqnarray*}
As we may get $\epsilon_3$ arbitrarily small, we have obtained the proposition.
\end{proof}
\begin{lemma}\label{G_2_lemma}
We have
\begin{gather*}
\sum\limits_{B \in \Delta}\dbeta(B)^2\diam(B) \lesssim 2^{-M\over 2}\cH^1 (\Gamma).
\end{gather*} 
\end{lemma}
\begin{proof}
For $B\in \Delta$ and $Q=Q(B)$, we will construct a weight $w_Q$  that satisfies (i), (ii) and (iii):  \\
\indent(i) $\int_Q w_Qd\ell\geq \diam(Q)$,\\
\indent(ii) for almost every $x_0\in \Gamma$, 
	$\sum\limits_{B \in \Delta} w_{Q(B)}(x_0) < C 2^{M\over 2}$,\\
\indent(iii) $\supp (w_Q) \subset Q$,\\
where $C$ is a constant which depends only on the 1-Ahlfors-regularity constant of $\Gamma$. 

%

We will construct $w_Q$ as a martingale.
We denote by $w_Q(Z):=\int_Z w_Qd\ell$.
Set
\begin{gather*}
w_Q(Q)=\diam(Q) .
\end{gather*}
Assume now that $w_Q(Q')$ is defined.  We define $w_Q(Q'^i)$ and $w_Q(R_{Q'})$,
where
\begin{gather*}
 Q'=(\cup  Q'^i) \cup R_{Q'},
\end{gather*} 
a decomposition as given by equation \eqref{05_11_05-7}.

Take
\begin{gather*}
w_Q(R_{Q'})=\frac{w_Q( Q')}{s'} \ell(R_{Q'}) 
\end{gather*}
and
\begin{gather*} 
w_Q( Q'^i)=\frac{w_Q(Q')}{s'}\diam(Q'^i),
\end{gather*}
where 
\begin{gather*}
s'=\ell(R_{Q'})+\sum_i \diam(Q'^i).
\end{gather*}
This will give us $w_Q$.
Note that $s'\lesssim \ell(\Gamma\cap Q')$.
Clearly (i) and (iii) are satisfied.
To see (ii):
\begin{eqnarray*}
\frac{w_Q( Q'^{i^*})}{\diam(Q'^{i^*}) }
&=&
\frac{w_Q( Q')}{s'}\\
&=&
\frac{w_Q( Q')}{\diam(Q') }
\frac{\diam(Q' )}{s'}\\
&=&
\frac{w_Q( Q')}{\diam(Q') }
\frac{\diam(Q') }
	{\ell(R_{Q'}) + 
		\sum\limits_{i} \diam(Q'^i) }\\
&\leq&
\frac{w_Q( Q')}{\diam(Q') }
\frac{1}
	{1+ c''\dbeta(B)}\\
\end{eqnarray*}

And so,
\begin{eqnarray*}
\frac{w_Q( Q'^{i^*})}{\diam(Q'^{i^*})} \leq 
 	q   \frac{w_Q( Q')}{\diam(Q')}
\end{eqnarray*}
With $q=\frac{1}
	{1+ c''2^{-M\over 2}}$\\
Now, suppose 
that  $x\in Q_N \subset ...\subset Q_1$.
we  get:
\begin{eqnarray*}
\frac{w_{Q_1}(Q_N)}{\diam(Q_N)} &\leq& 
  q\frac{w_{Q_1}(Q_{N-1})}{\diam(Q_{N-1})} \\
  &\leq&...\\
  &\leq&
  q^{N-1}\frac{w_{Q_1}(Q_{1})}{\diam(Q_1)}=q^{N-1}.
\end{eqnarray*}
 Hence,  we have $w_{Q_1}(x) \lesssim q^{N-1}$.
This will give us (ii) as a sum of a geometric series'  since
\begin{gather*}
\sum q^n = {1 \over 1-q}\lesssim {1\over 2^{-M\over 2}}=2^{M\over 2}.
\end{gather*}

Now,
\begin{eqnarray*}
\sum\limits_{B \in \Delta}\dbeta(B)^2\diam(B) 
&\lesssim& 
2^{-M}\sum\limits_{B \in \Delta}\diam(B)\\
&\lesssim&
2^{-M}\sum\limits_{B \in \Delta}\int w_{Q(B)}(x)d\ell(x)\\
&=&
2^{-M}\int \sum\limits_{B \in \Delta} w_{Q(B)}(x)d\ell(x)\\
&\lesssim&
2^{-M}\int  2^{M\over 2} d\ell(x)\\
&\lesssim&
2^{-M\over 2}\cH^1(\Gamma).
\end{eqnarray*}
\end{proof}
\begin{rem}
By taking an increasing sequence of $\Delta\to \Delta^{p_1,p_2}_M$ we get that  Lemma \ref{G_2_lemma} holds with $\Delta$ is replaced by $\Delta^{p_1,p_2}_M$.
\end{rem}

We now get Theorem \ref{new-thm-2} since 
\begin{eqnarray*}
&&\sum\limits_{B \in \G_2}\beta^2_2(B)^2\diam(B)\\
&\lesssim& 
\sum\limits_{B \in \G_2}\dbeta(B)^2\diam(B)\\ 
&=&
\sum\limits_{M=1}^\infty 
	\sum\limits_{B \in \G_2 \atop 2^{-M} \leq \frac12 \dbeta(B)^2  < 2^{-M+1}}
		\dbeta(B)^2\diam(B) \\
&\lesssim& 
\sum\limits_{M=1}^\infty 
	M2^{-M\over 2}\cH^1(\Gamma)\\
&\lesssim& 
\cH^1(\Gamma)\,.
\end{eqnarray*}

%% file: local_var/local_var.tex
%
%
In this section we give the needed modifications to obtain Theorems 
 \ref{new-thm-1-var} and \ref{new-thm-2-var}.

Consider a ball $\ball(z,R)$ where $R>0$ and $z\in \Gamma$.
Let $\{\Gamma_i\}$ be the  connected components of $\Gamma \cap \ball(z,10R)$ which intersect $\ball(z,R)$. If there is only one such component then  Theorems \ref{new-thm-1} and  \ref{new-thm-2} give Theorems \ref{new-thm-1-var} and \ref{new-thm-2-var}.
Otherwise, all  components $\Gamma_i$ must have diameter at least $9R$, and so by 1-Ahlfors-regularity there are at most $P$ of them, where $P$ depends only on the 1-Ahlfors-regularity constant of $\Gamma$.
Parameterize each $\Gamma_i$  by $\gamma_i$, as assured by Lemma \ref{parameterization}.

{\it Informally} speaking, the proofs we have of Theorems \ref{new-thm-1} and \ref{new-thm-2} now work word for word, since they only depend on the existence of a parameterization for each connected component.  
Rather than checking this, we use the following trick.
 
One may simply connect the end of $\gamma_i$ to the beginning of $\gamma_{i+1}$ with an arc-length parameterization.  The total added length will be at most $20PR$.  Call this new path $\gamma$, and its image $\tilde{\Gamma}$.  One may apply Theorems \ref{new-thm-1} and \ref{new-thm-2} to get the desired results now.

This completes the proof of Theorems \ref{new-thm-1-var} and \ref{new-thm-2-var}.

%% file: appendix/parameterization.tex
\subsection{Proof of Lemma \ref{parameterization}}
We assume  $\Gamma \subset \M$ be a compact connected  set of finite $\cH^1$ measure.

Using the 
Kuratowski embedding theorem (see \cite{Heinonen_embedding_lec}), 
we have  an isometric embedding $f:\Gamma\to \ell_\infty(\Gamma)$.
Let $\Gamma'=f(\Gamma)$.  

The following two lemmas have proofs identical to what appears in \cite{my-thesis-as-paper}.
\begin{lemma}
Let $C_1,C_2 >0$ be given.
Given a compact  set $\Gamma' \subset \ell_\infty(\Gamma)$ the set
$E:=\{x\in \ell_\infty(\Gamma):x=tx_1 + (1-t)x_2, x_i \in \Gamma', -C_1\leq t \leq C_2\}$
is compact.
\end{lemma}
\begin{proof}
Suppose $\{x^i\} \subset E$ is a sequence.  
We can write $x^i=t^ix^i_1 + (1-t^i)x^i_2$ as in the definition of $E$.
By the compactness of $\Gamma'$ we have 
$i_k \st \\
x^{i_k}_1 \to x_1$. By compactness of $\Gamma'$ again, $x^{i_{k_j}}_2 \to x_2$. 
By  compactness of  $[-C_1,C_2]$ we have $t^{i_{k_{j_l}}} \to t$. 
$x_1,x_2 \in \Gamma', t \in [-C_1,C_2]$.  Hence $x^{i_{k_{j_l}}}\to tx_1+ (1-t)x_2 \in E$.
\end{proof}

\smallskip
\noindent
\begin{lemma}
Let $\Gamma' \subset \ell_\infty(\Gamma)$ be a compact connected  set of finite length.  
Then we have a Lipschitz function $\gamma:[0,1] \to \ell_\infty(\Gamma)$ such that $Image (\gamma)=\Gamma'$ 
and $\norm{\gamma}_{Lip} \leq 32\cH^1 (\Gamma')$
\end{lemma}
\begin{proof}
We use a  well known result from graph theory (which we call (*)):\\
If $G$ is a connected graph with finitely many edges, then there is a path that traverses each 
edge of $G$ exactly twice (once in each direction).
This result is easily seen by induction on the number of edges.

For $n\geq0$, let $X_n=X_n^{\Gamma'}$ (i.e. take $X_n \subset \Gamma'$ a $2^{-n} - net$).
Note that since $\Gamma'$ is compact, each $X_n$ is finite.
We want to get a connected set $E_n$.
We do this by adding line segments inductively.
Set $E_n^0=X_n$.  
We get  $E_n^{i+1}$ from $E_n^i$ by adding a line segment between points 
$x_1,x_2 \in X_n$ such that $\dist(x_1,x_2)<2^{-n+3}$ and they are not yet in the same 
connected component of $E_n^i$.  
If there are no two such points we stop and call the resulting set $E_n$.
Let $G_n$ be the obvious abstract graph associated to $E_n$. 
If $G_n$ is not connected then ${\rm Vertex} (G_n)=A \cup B$ with $\dist(A,B) \geq 2^{-n+2}$ and  
$A$ separated from $B$.
By the construction of $E_n$ and $X_n$ we have that $\dist(\N_{2^{-n}}(A),\N_{2^{-n}}(B)) \geq 2^{-n}$ and $\Gamma' \subset \N_{2^{-n}}(A)\cup \N_{2^{-n}}(B)$.  This is a contradiction to  $\Gamma'$ being connected.  
Hence $G_n$ is connected.\\
Note that $\cH^1(E_n) \leq \sharp(X_n)2^{-n+3} \leq 16 \cH^1(\Gamma')$, where the final inequality follows from the fact that the balls $\{B(x,2^{-n-1}):x \in X_n\}$ are disjoint.\\
We can thus parameterize $E_n$ by a Lipschitz curve  $\gamma_n:[0,1] \to \ell_\infty(\Gamma)$. 
The image of this parameterization is in $E$ as defined in the previous lemma.  
By Arzela-Ascoli we have a subsequence converging to $\gamma'$.
We have that $Image(\gamma')=\Gamma'$ by say 
\begin{gather*}
\supl_{x\in E_n} \dist(x,\Gamma')+ \supl_{y\in \Gamma'} \dist(E_n,y) \leq 4\cdot 2^{-n} +   2^{-n}
= 5\cdot 2^{-n}
\end{gather*}
 and a triangle inequality. 
\end{proof}

Now,  Consider the mapping $\gamma''=f^{-1}\gamma'$.  The map $\gamma''$ gives the first part of Lemma \ref{parameterization} with $\T$ replaced by $[0,1]$. To correct this one simply defines 
$\gamma(t)=\gamma''(2t)$ for $0\leq t\leq \half$ and 
$\gamma(t)=\gamma''(1-(2t-1))$ for $\half \leq t\leq1 $. The map $\gamma$ has $\T$ as its domain and $\Gamma$ as its image.

Assume now that $\Gamma$ is also 1-Ahlfors-regular with constant $C$.
Then in the proof above, $E_n$ is also 1-Ahlfors-regular.  
Hence
$$\frac{R}{C'}\leq\cH^1(\gamma_n^{-1}(\ball(x,R)))\leq C'R \quad
\forall x\in E_n,\ 0<R\leq\diam(E_n)$$
by the result (*). Given $R$, one may choose $n$ large enough so that this implies the second part of  
Lemma \ref{parameterization}.   

This completes the proof of Lemma \ref{parameterization}.
%